\documentclass[12pt]{article}
\usepackage{amsfonts}

\usepackage{amssymb}
\usepackage{amsmath}
\usepackage{endnotes}
\usepackage[doublespacing]{setspace}

\newenvironment{proof}[1][Proof]{\noindent\textbf{#1.} }{\ \rule{0.5em}{0.5em}}
\input{tcilatex}
\let\footnote=\endnote

\begin{document}

\title{Wald tests when restrictions are locally singular\thanks{%
This work was supported by the Willam Dow Chair in Political Economy (McGill
University), the Bank of Canada Research Fellowship, The Toulouse School of
Economics Pierre-de-Fermat Chair of Excellence, A Guggenheim Fellowship,
Conrad-Adenauer Fellowship from Alexander-von-Humboldt Foundation, the
Canadian Network of Centres of Excellence program on Mathematics of
Information Technology and Complex Systems, the Natural Sciences and
Engineering Research Council of Canada, the Social Sciences and Humanities
Research Council of Canada and the Fonds de recherche sur la soci\'{e}t\'{e}
et la culture (Qu\'{e}bec). The authors also thank the research centres
CIREQ and CIRANO for providing support and meeting space for the joint work.
We thank Purevdorj Tuvaandorj for very useful comments.}}
\author{Jean-Marie Dufour\thanks{%
William Dow Professor of Economics, McGill University, Centre
interuniversitaire de recherche en analyse des organisations (CIRANO) and
Centre interuniversitaire de recherche en \'{e}conomie quatative (CIREQ).}
\and Eric Renault \thanks{%
Brown University} \and Victoria Zinde-Walsh\thanks{%
McGill University and CIREQ} \\
}
\maketitle
\date{}

\textbf{Key words}:\vspace{0.25pt}nonlinear restriction; deficient rank;
singular covariance matrix; Wald test\textbf{;} nonstandard\textbf{\ }%
asymptotic theory; bound. \textbf{Journal of Economic Literature\
classification:} C3.

\bigskip \pagebreak

\section{\protect\bigskip Introduction}

Tests based on asymptotic distributions typically require regularity
assumptions in order to be able to obtain critical values. This is the case,
in particular, for Wald-type statistics based on asymptotically normal
estimators. Wald-type tests are especially convenient because they allow one
to test a wide array of linear and nonlinear restrictions from a single
unrestricted estimator. We focus here on the problem of implementing
Wald-type tests for nonlinear restrictions.

The use of the Wald statistic has been criticized because of finite sample
non-invariance (Gregory and Veall (1985), Breusch and Shmidt (1988),
Phillips-Park(1988), Dagenais-Dufour(1991)) and lack of robustness to
identification failure (Dufour(1997,2003)). We focus here on situations
where the parameter tested is typically identified under the null
hypothesis, but usual rank conditions on the Jacobian matrix may fail
asymptotically.

Under regularity conditions, the standard asymptotic distribution of the
test statistic is chi-square with degrees-of-freedom equal to the number of
restrictions. The regularity conditions involve the assumption that the
restrictions are differentiable with respect to the parameters considered,
with a derivative matrix which has full column rank in an open neighborhood
of the true value of the parameter vector. There are many problems, however,
for which this regularity condition is violated. These include, among others:

\begin{enumerate}
\item hypothesis tests on bilinear and multilinear forms of model
coefficients in Gourieroux-Monfort-Renault (1988);

\item testing whether the matrix of polynomials or multilinear forms in
model coefficients has than full rank or, equivalently, whether the
determinant of the matrix is zero in Gourieroux-Monfort-Renault (1993);

\item tests of Granger noncausality in VARMA models in
Boudjellaba-Dufour-Roy (1992,1994);

\item tests of noncausality at various horizons in Dufour-Renault (1998),
Dufour-Pelletier-Renault (2005);

\item tests for common factors in ARMA models in Gourieroux-Monfort-Renault
(1989), Galbraith-ZindeWalsh (1997);

\item test of volatility and covolatility\ in Gouri\'{e}roux and Jasiak
(2013).
\end{enumerate}

A common feature of the above problems is the fact that the estimated
asymptotic covariance matrix of the relevant nonlinear functions of
coefficient estimates converges to a singular matrix on a subset of the null
hypothesis -- so that the usual regularity condition fails -- but is
non-singular (with probability one) in finite samples. The estimated
covariance matrix used by the Wald-type statistic is a consistent estimator
of the asymptotic covariance matrix of the corresponding nonlinear form in
parameter estimates, but the rank of the estimated covariance matrix does
not consistently estimate the rank of the asymptotic covariance matrix
(because the rank is not a continuous function). It is important to note
here that this is \emph{not an identification problem}, so that standard
criticisms of Wald-type methods in the presence of identification problems
(see Dufour (1997,2003)) do not apply in this case.

If the covariance matrix estimator can be modified so that it remains
consistent and its rank converges to the appropriate asymptotic rank, then
the asymptotic distribution of the modified Wald-type statistic (based on a
generalized inverse of the covariance matrix) remains chi-square although
with a reduced degrees-of-freedom number; see Andrews (1987). For example,
Lutkepohl-Burda (1997) proposed such methods based on reducing the rank of
the estimated covariance matrix by either using a form of randomization or
setting \textquotedblleft small eigenvalues\textquotedblright\ to zero. Such
methods, however, effectively modify the test statistic and involve
arbitrary truncation parameter for which no practical guidelines are
available: in finite samples, the test statistic can become as small as one
wishes leading to largely arbitrary results and unlimited power reductions.

Interestingly, except for a bound given by Sargan(1980) in a special case,
the asymptotic distribution of Wald-type statistics in non-regular cases has
not been studied. In this paper, we undertake this task and propose
solutions to the problem that do not require modifying the test statistic.
More specifically, the contributions of the paper can be summarized as
follows.

\emph{First}, we provide examples showing that Wald statistics in such
non-regular cases can have several asymptotic distributions. We also show
that usual critical values based on a chi-square distribution (with
degrees-of-freedom equal to the number of constraints) can both lead to
under-rejections and over-rejections depending on the form of the function
studied. Indeed, the Wald statistic may diverge under null hypothesis, so
that arbitrary size distortions may occur.

\emph{Second}, we study the asymptotic distribution of Wald-type statistics
in non-regular cases. Surprisingly, the asymptotic behavior of the Wald
statistic has not been generally studied for full classes of restrictions;
here we consider the class of polynomial restrictions. We show that the Wald
statistic either has a non-degenerate asymptotic distribution even when the
estimated covariance converges to a singular matrix, or diverges to
infinity. We provide conditions for convergence and a general
characterization of this distribution. We find that the test can have
several different asymptotic distributions under the null hypothesis --
depending on the degree of singularity as well as various nuisance
parameters -- which may be non-chi-square distributions.

\emph{Third}, we provide bounds on the asymptotic distribution (when it
exists), which turn out to be to be proportional to a chi-square
distribution where the proportionality constant depends on the degree of
singularity of the function considered. In several cases of interest, this
bound yields an easily available conservative critical value. Even when the
limit distribution is non-pivotal it is sometime possible to provide pivotal
bounds that would yield conservative critical values.

\emph{Fourth}, we propose an adaptive consistent strategy for determining
whether the asymptotic distribution exists and which form it takes; this
approach also permits to determine what kind of bound is valid.

The framework considered and the test statistics are defined in Section 2. A
number of examples are presented in Section 3; they illustrate the
properties of the Wald test in singular cases. In Section 4 we discuss some
general algebraic and analytic features of matrices of polynomials and
quadratic forms and derive the asymptotic distribution of the Wald
statistic. Bounds are derived in Section 5. An adaptive strategy for
determining the asymptotic distribution and the bounds is developed in
Section 6. Proofs are presented in the Appendix.

\section{Framework}

We consider testing $q$ restrictions in a situation where an asymptotically
non-singular estimator $\hat{\theta}_{T}$ is available for a $p\times 1$
parameter of interest $\bar{\theta}$ that satisfies the restrictions; $q\leq
p.$

Assumption 2.1.\textit{\ The function }$g\left( \theta \right) =[g_{1}\left(
\theta \right) ,\ldots ,\,g_{q}\left( \theta \right) ]^{\prime }$\textit{\
is a continuously differentiable function from }$\Theta $\textit{\ to }$R^{q}
$\textit{, where }$\Theta $\textit{\ is an open subset of }$R^{p}$\textit{\
and }$q\leq p.$

Assumption 2.1a. \textit{The function }$g\left( \theta \right) =[g_{1}\left(
\theta \right) ,\ldots ,\,g_{q}\left( \theta \right) ]^{\prime }$\textit{\
is such that each }$g_{i}\left( \theta \right) $\textit{\ is a polynomial of
order }$m$\textit{\ in the components of }$\theta ,$\textit{\ i.e.}%
\begin{eqnarray}
g_{i}\left( \theta \right)  &=&\underset{k=0}{\overset{m}{\sum }}%
g_{ik}(\theta )\,,  \label{eq: gik(theta)} \\
g_{ik}(\theta ) &=&\underset{i_{1}+\cdots +i_{p}=k}{\sum }%
A_{ik}(i_{1},\ldots ,\,i_{p})\,\theta _{1}^{i_{1}}\cdots \,\theta
_{p}^{i_{p}}\,,\;k=0,1,\ldots ,\,m,\;i=1,\ldots ,\,q,
\label{eq: gik(theta) polynomial}
\end{eqnarray}%
\textit{where }$g_{ik}(\theta )$\textit{\ represents a homogeneous
polynomial of order }$k$\textit{, each coefficient }$A_{ik}(i_{1},\ldots
,\,i_{p})$\textit{\ is a constant, and }$m$\textit{\ is the maximal order of
a polynomial in }$g\left( \theta \right) $.

Assumption 2.2. \textit{We assume that some }$\bar{\theta}$\textit{\
satisfies a null hypothesis of the form: }%
\begin{equation}
H_{0}:g\left( \theta \right) =0\,.  \label{hyp}
\end{equation}

Assumption 2.3. \textit{Assume that }$\{\hat{\theta}_{T}:$\textit{\ }$T\geq
T_{0}\}$\textit{\ is a sequence of }$p\times 1$\textit{\ random vectors such
that for some positive definite matrix }$V$\textit{\ and a scalar rate
sequence }$\lambda _{T}\rightarrow \infty $\textit{\ as }$T\rightarrow
\infty $\textit{\ convergence in probability holds: }%
\begin{equation}
\lambda _{T}V^{-\frac{1}{2}}\left( \hat{\theta}_{T}-\bar{\theta}\right)
\rightarrow _{p}Z,  \label{eq: Asymptotic distribution}
\end{equation}%
\textit{where }$Z$\textit{\ is a random }$p\times 1$\textit{\ vector with a
known absolutely continuous probability distribution, }$Q\left( \bar{\theta}%
\right) $\textit{\ on }$R^{p}.$

Assumption 2.3a.\textit{\ In addition to Assumption 2.3 }$\lambda _{T}=T^{%
\frac{1}{2}},$\textit{\ }$Z$\textit{\ is a Gaussian random vector.}

Assumption 2.4. $\{\hat{V}_{T}:$\textit{\ }$T\geq T_{0}\}$\textit{\ is a
sequence of }$p\times p$\textit{\ random matrices such that }$P[rank(\hat{V}%
_{T})=p]=1,$\textit{\ for all }$T,$\textit{\ and}%
\begin{equation}
\underset{T\rightarrow \infty }{\mathrm{plim}}\,\hat{V}_{T}=V
\label{eq: Consistent V estimator}
\end{equation}%
\textit{where the probability that }$\hat{V}_{T}$\textit{\ be positive
definite is one for }$T\geq T_{0}$\textit{\ }$($\textit{for some }$T_{0}>0)$.

We define the Wald test statistic:

\begin{equation}
W_{T}=\lambda _{T}^{2}g^{\prime }(\hat{\theta}_{T})\left[ \frac{\partial g}{%
\partial \theta ^{\prime }}(\hat{\theta}_{T})\hat{V}_{T}\frac{\partial
g^{\prime }}{\partial \theta }(\hat{\theta}_{T})\right] ^{-1}g(\hat{\theta}%
_{T}),  \label{waldgen}
\end{equation}%
when $\lambda _{T}^{2}=T,$ this is%
\begin{equation}
W_{T}=Tg^{\prime }(\hat{\theta}_{T})\left[ \frac{\partial g}{\partial \theta
^{\prime }}(\hat{\theta}_{T})\hat{V}_{T}\frac{\partial g^{\prime }}{\partial
\theta }(\hat{\theta}_{T})\right] ^{-1}g(\hat{\theta}_{T}).
\label{waldstand}
\end{equation}%
If the distribution $Q(\bar{\theta})$\ has a finite variance, we can assume
without loss of generality that its variance is the identity matrix.

However, when the rate of convergence $\lambda _{T}$\ is not the standard $%
T^{1/2}$, a factor $\lambda _{T}^{2}$\ shows up instead of $T$.

The statistic $W_{T}$ is not well defined when the estimator $\hat{\theta}%
_{T}$ falls into the set of singularity points at which $\frac{\partial g}{%
\partial \theta ^{\prime }}(\hat{\theta}_{T})\hat{V}_{T}\frac{\partial
g^{\prime }}{\partial \theta }(\hat{\theta}_{T})$ is non-invertible (of rank
less than $q)$. Andrews (1987) studied the case where $\left[ \frac{\partial
g}{\partial \theta ^{\prime }}(\hat{\theta}_{T})\hat{V}_{T}\frac{\partial
g^{\prime }}{\partial \theta }(\hat{\theta}_{T})\right] ^{-1}$ is replaced
by a generalized inverse (\emph{e.g.}, the Moore-Penrose inverse) and gave
conditions under which the asymptotic distribution is chi-square. The main
result there is that the asymptotic distribution of $W_{T}$ under $H_{0}$ is
chi-square $\chi ^{2}(r_{0})$ with $r_{0}=$\textrm{$rank$}$[\frac{\partial g%
}{\partial \theta ^{\prime }}(\theta )]$ when \textrm{$rank$}$[\frac{%
\partial g}{\partial \theta ^{\prime }}(\hat{\theta}_{T})]$ converges to $%
r_{0}$ under $H_{0}.$ This will be the case in particular when $\frac{%
\partial g}{\partial \theta ^{\prime }}(\theta )$ has rank $r_{0}$ in some
open neighborhood of $\bar{\theta}.$

Here we study situations where the matrix $\frac{\partial g}{\partial \theta
^{\prime }}(\hat{\theta}_{T})\hat{V}_{T}\frac{\partial g^{\prime }}{\partial
\theta }(\hat{\theta}_{T})$ is non-singular in finite samples (with
probability 1) but may converge to a singular matrix. Under Assumption 2.4
this non-singularity is equivalent to the matrix $G\left( \theta \right) =%
\frac{\partial g}{\partial \theta ^{\prime }}(\theta )$ having full rank
almost everywhere.

Assumption 2.5. \textit{The matrix }$G\left( \theta \right) $\textit{\ has
full row rank for almost all }$\theta .$

\section{Examples and counter-examples}

Before we move to study the asymptotic distribution of $W_{T}=W_{T}(\hat{%
\theta}_{T},\,\hat{V}_{T})$ in general terms we provide examples which show
that, indeed, the asymptotic distribution of $W_{T}$ is not regular. In
particular, our examples illustrate non-invariance of the asymptotic
distribution of the statistic to the form of the restriction and dependence
(discontinuous) of the asymptotic distribution on the parameter value, $\bar{%
\theta};$ we also show that the asymptotic distribution may have either
thinner or thicker tails than the standard $\chi _{q}^{2}$ distribution and
can even diverge to infinity under the null.

To streamline exposition of the examples we assume that $V=I.$

The following example illustrates lack of invariance of the asymptotic
distribution.

\textbf{Example 3.1. }\textit{Consider two equivalent forms for the null, }$%
g\left( \theta \right) =0,$\textit{\ one is (i) }$\theta =\allowbreak 0,$%
\textit{\ the other (ii) }$\theta ^{2}=0.$\textit{\ Of course, the
asymptotic distribution for the Wald test statistic in the case (i) under
Assumption 2.3a is }$\chi _{1}^{2}.$\textit{\ By contrast, for (ii) the
value of }$W_{T}=T\frac{\hat{\theta}_{T}^{4}}{4\hat{\theta}_{T}^{2}};$%
\textit{\ the limit distribution then is }$\frac{1}{4}\chi _{1}^{2}.$

Below for the multivariate $\hat{\theta}_{T}$ we suppress dependence of the
components, $\hat{\theta}_{T1},...,\hat{\theta}_{Tp},$ on $T.$

The next example is the one given by Andrews (1987); we develop it to
illustrate both the fact that the distribution depends on $\bar{\theta},$
and also that despite the distribution not being pivotal, the usual $\chi
_{1}^{2}$ distribution provides here a pivotal upper bound.

\textbf{Example 3.2. }\textit{Consider the restriction given by }$g(\theta
)=\theta _{1}\theta _{2}.$\textit{\ In this case, }$G\left( \theta \right)
=[\theta _{2},\,\theta _{1}],$\textit{\ and the Wald statistic for testing }$%
H_{0}:\theta _{1}\theta _{2}=0$\textit{\ takes the form: }%
\begin{equation*}
W_{T}=T\frac{\hat{\theta}_{1}^{2}\hat{\theta}_{2}^{2}}{\hat{\theta}_{1}^{2}+%
\hat{\theta}_{2}^{2}}\,.
\end{equation*}%
\textit{If either }$\theta _{1}$\textit{\ or }$\theta _{2}$\textit{\ is
non-zero, under }$H_{0}$\textit{\ the limiting distribution is }$\chi
_{1}^{2}.$\textit{\ If, however, }$\theta _{1}=\theta _{2}=0,$\textit{\ we
have:}%
\begin{equation}
W_{T}\overset{p}{\underset{T\rightarrow \infty }{\longrightarrow }}\frac{%
Z_{1}^{2}Z_{2}^{2}}{Z_{1}^{2}+Z_{2}^{2}}.  \label{andrews}
\end{equation}

\textit{Writing this expression as }$Z_{2}^{2}-\frac{Z_{2}^{4}}{%
Z_{1}^{2}+Z_{2}^{2}}$\textit{\ we see that the limit distribution in this
case under Assumption 2.3a is strictly below }$\chi _{1}^{2},$\textit{\ thus
it is not pivotal. However, }$\chi _{1}^{2}$\textit{\ provides a
conservative bound.}

\textit{A more precise bound can be obtained. Write the vector }$\left(
Z_{1},Z_{2}\right) $\textit{\ in polar coordinates: }$\left( r\sin \phi
,r\cos \phi \right) ,$\textit{\ with }$r^{2}=Z_{1}^{2}+Z_{2}^{2},r\geq 0$ 
\textit{and} $\phi =\arcsin \frac{Z_{1}}{r}.$ \textit{Then the limit ratio
in }$\left( \ref{andrews}\right) $\textit{\ becomes}%
\begin{equation*}
\frac{1}{4}r^{2}\left( \sin 2\phi \right) ^{2}.
\end{equation*}%
\textit{Thus the distribution of }$\frac{1}{4}r^{2}$\textit{\ provides an
upper bound on the limit distribution of }$W_{T}$\textit{\ under the most
general assumptions.}

\textit{If the distribution of the vector }$Z$\textit{\ is spherical (that
is depends on }$r$ only)\textit{, then the distribution of }$\phi $\textit{\
is uniform and independent of }$r$\textit{; it follows that }$r\sin 2\phi $%
\textit{\ has then the same distribution as }$r\sin \phi .$\textit{\ Indeed,
conditionally on }$r$\textit{\ (denoting by }$F_{\cdot |\cdot }\left( \cdot
\right) $\textit{\ the conditional distribution) }%
\begin{eqnarray*}
F_{\sin 2\phi |r}(\frac{\alpha }{r}) &=&F_{2\phi |r}\left( \arcsin \frac{%
\alpha }{r}\right) =2F_{\phi |r}\left( \frac{1}{2}\arcsin \frac{\alpha }{r}%
\right)  \\
&=&2\int_{0}^{\frac{1}{2}\arcsin \left( \alpha /r\right) }I\left( 0\leq \phi
\leq 2\pi \right) \frac{1}{2\pi }d\phi =F_{\phi |r}\left( \arcsin \frac{%
\alpha }{r}\right) .
\end{eqnarray*}%
\textit{Then the limit of }$W_{T}$\textit{\ is given by the distribution of }%
$\frac{1}{4}Z_{1}^{2}$\textit{\ (the same as }$\frac{1}{4}Z_{2}^{2}).$

\textit{Under normality this is distributed as }$\frac{1}{4}\chi _{1}^{2}.$%
\textit{\ If the distribution of }$Z$\textit{\ is such that each marginal is
normal but the joint is not, then }$\frac{1}{4}\chi _{1}^{2}$\textit{\
provides an upper bound but not necessarily the distribution.}

The limit $\frac{1}{4}\chi _{1}^{2}$ distribution under normality was
obtained by Glonek (1993) who also demonstrated that this asymptotic
distribution does not depend on the covariance matrix $V.$ Thus the limit
distribution for test of this hypothesis for a normal $Z$ is either $\chi
_{1}^{2}$ or $\frac{1}{4}\chi _{1}^{2}$, therefore is not pivotal. However, $%
\chi _{1}^{2}$ provides a conservative bound, so that there is a pivotal
upper bound.

In the above examples, standard critical values are conservative in
non-regular cases. So here if we do not know whether we are in a regular
case or not, usual critical values are the appropriate ones: the test never
over-rejects (asymptotically) under the null hypothesis when using critical
values entailed by usual regularity assumptions.

However, it is also possible that the standard limit distribution does not
hold in any part of the parameter space and using the corresponding critical
values may lead to a severely oversized test.

\textbf{Example 3.3. }\textit{Suppose that }$g\left( \theta \right) =\theta
_{1}^{2}+...+\theta _{p}^{2};$\textit{\ then }$G\left( \theta \right) =\left[
2\theta _{1},....,2\theta _{p}\right] $\textit{\ and }%
\begin{equation*}
W_{T}=T\frac{\left( \Sigma _{i=1}^{p}\hat{\theta}_{i}^{2}\right) ^{2}}{%
4\Sigma _{i=1}^{p}\hat{\theta}_{i}^{2}}.
\end{equation*}

\textit{Then the limit distribution is that of }$\frac{1}{4}\left\Vert
Z\right\Vert ^{\frac{1}{2}};$\textit{\ under normality this is }$\frac{1}{4}%
\chi _{p}^{2};$\textit{\ it is a pivotal distribution even though
non-standard. If }$p$\textit{\ is large enough, the }$\chi _{1}^{2}$\textit{%
\ will not provide an upper bound.}

In the case of more than one restriction in addition to all the non-standard
features that can arise for a single restriction it is also possible that
the test statistic diverges even under $H_{0}.$

\textbf{Example 3.4. }\textit{Suppose that }$q=p=2$\textit{\ and }$g\left(
\theta \right) =\left[ \theta _{1}^{2}:\theta _{1}\theta _{2}^{2}\right]
^{\prime }.$\textit{\ Then }%
\begin{equation*}
G\left( \theta \right) =\left[ 
\begin{array}{cc}
2\theta _{1} & 0 \\ 
\theta _{2}^{2} & 2\theta _{1}\theta _{2}%
\end{array}%
\right] ;
\end{equation*}%
\textit{it follows that }%
\begin{equation*}
W_{T}=T\frac{4\hat{\theta}_{1}^{2}+\hat{\theta}_{2}^{2}}{16}.
\end{equation*}%
\textit{Then if (i) }$\bar{\theta}_{1}=\bar{\theta}_{2}=0$\textit{\ the
asymptotic distribution is }$\frac{1}{4}Z_{1}^{2}+\frac{1}{16}Z_{2}^{2}$%
\textit{\ and thus under normality is a linear combination of two
independent }$\chi _{1}^{2}$\textit{\ and is bounded by }$\frac{1}{4}\chi
_{2}^{2}.$\textit{\ However, if (ii) }$\bar{\theta}_{1}=0,$\textit{\ but }$%
\bar{\theta}_{2}\neq 0$\textit{\ the null still holds, but as }$T\rightarrow
\infty $\textit{\ the Wald statistic diverges to }$+\infty .$

The examples show that even for the simplest restrictions the limit
distribution of the Wald statistic may be quite complex and far from
standard. A number of applications require the Wald test of polynomial
restriction functions where singularity could not be excluded and thus the
non-standard features illustrated by the simple examples above may be
present.

Several applications involve test of one restriction, such as tests of
determinants and other polynomial functions in coefficients in Gourieroux,
Monfort, Renault (1988, 1993), Galbraith and Zinde-Walsh (1992), Gourieroux
and Jasiak (2013). In tests of Granger noncausality in VARMA models by
Boudjellaba, Dufour and Roy (1992,1994) several polynomial restrictions need
to hold under the null, similarly in testing noncausality at various
horizons in Dufour, Renault (1998) and Dufour, Pelletier and Renault (2005).

\section{Limit distribution of the Wald statistic}

The asymptotic behavior of the Wald statistic has not been generally
examined in the literature for full functional classes of nonlinear
restrictions. Here we provide a characterization of the asymptotic
distribution for restrictions given by polynomial functions. We shall work
under the Assumptions 2.1a, 2.2,2.3, 2.4 and 2.5.

Two approaches are possible. The Wald statistic can be represented as a
ratio of two polynomial functions in random variables; such a representation
implicitely incorporates the information in the polynomial restrictions.
Another approach is based on an explicit analysis of the restrictions and
represents the limit distribution in a quadratic form; this representation
permits simple derivation of conservative bounds. In this paper we focus on
the second representation.

The first subsection gives a few general results about matrices of
polynomials; the second applies them to matrices related to the Jacobian
matrix of the restrictions under test. The third subsection provides the
limit distribution for the Wald statistic for polynomial restrictions; this
distribution is in general not pivotal and depends on $\bar{\theta}$.

\subsection{Matrices of polynomials}

A polynomial function is either the zero polynomial, when it is identically
zero (the coefficient on every monomial term is zero), or it is non-zero
a.e..

Consider a $q\times p$\ matrix $G(y)$ of polynomials of variable $y\in 
\mathbb{R}^{p}.$ When $q=p,$ we will say that the matrix $G(y)$\ is
non-singular if its determinant is a non-zero polynomial. More generally, we
will define the rank of the $q\times p$\ matrix $G(y)$ as the largest
dimension of a square non-singular submatrix. This section considers\ $%
q\times p$\ matrices $G(y),$\ $q\leq p,$ of full row rank $q$ (Assumption
2.5).

We first note that, for any square $q\times q$ non-singular matrix $S$, $%
SG(y)$ is also a matrix of polynomials of rank $q$: if $\breve{G}(y)$\ is a $%
q\times q$\ submatrix \ of $G(y)$\ with determinant $\det \left( \breve{G}%
(y)\right) $\ that is a non-zero polynomial, it is also true for the
submatrix $S\breve{G}(y)$\ of $SG(y)$.

Consider a polynomial $h(y)=\Sigma _{k=0}^{n}h_{k}(y)$ with homogeneous
polynomial terms of order $k:$ 
\begin{equation}
h_{k}(y)=\Sigma
_{i_{1}+...+i_{p}=k}h_{k}(i_{1},...,i_{p})y_{1}^{i_{i}}...y_{p}^{i_{p}}.
\label{polyn}
\end{equation}

Denote by $\bar{k}_{h}$ the lowest order of homogeneous polynomial entering
into polynomial $h(y):$%
\begin{equation}
\bar{k}_{h}=\underset{0\leq k\leq n}{\min }\left\{ k:h_{k}\left(
i_{1},...,i_{p}\right) \neq 0\text{ for some }i_{1}+...+i_{p}=k\right\} .
\label{k-bar}
\end{equation}

Note that

\begin{equation}
\lambda ^{\bar{k}_{h}}h(y/\lambda )=h_{\bar{k}_{h}}(y)+\Sigma \lambda
^{r_{l}}r_{l}\left( y\right) ,  \label{homog}
\end{equation}%
with all $r_{l}<0$ and $r_{l}(y)$ polynomial with $\bar{k}_{r_{l}}>\bar{k}%
_{h}.$

Consider all possible $\tilde{G}(y)_{l},$ with $\tilde{G}(y)_{l}$ a $q\times
q$ submatrix of $G\left( y\right) ;$ $l=1,...L$ with $L=\frac{p!}{q!(p-q)!}.$

Define 
\begin{equation}
\bar{\alpha}=\underset{l}{\min }(\bar{k}_{\det (\tilde{G}(y)_{l})})
\label{a-bar}
\end{equation}%
with the convention $\bar{k}_{\det (\tilde{G}(y)_{l})}=+\infty $ if $\det
\left( \tilde{G}(y)_{l}\right) $ is the zero polynomial.

Note that for some $\tilde{G}_{l}$ strict inequality $\bar{k}_{\det (\tilde{G%
}(y)_{l})}>\bar{\alpha}$ may hold as shown in the example below.

\textbf{Example 4.1. }$G\left( y\right) =\frac{\partial g}{\partial
y^{\prime }}$\textit{\ for }$g\left( y\right) =\left(
y_{1}^{2}+y_{3}^{3},y_{2}^{2}+y_{4}^{3},y_{1}^{2}+y_{2}^{2}\right) ^{\prime
};$\textit{\ then}

\begin{equation*}
\mathit{G}\left( y\right) \mathit{=}\left[ 
\begin{array}{cccc}
2y_{1} & 0 & 3y_{3}^{2} & 0 \\ 
0 & 2y_{2} & 0 & 3y_{4}^{2} \\ 
2y_{1} & 2y_{2} & 0 & 0%
\end{array}%
\right] \mathit{.}
\end{equation*}

\textit{We have four possible }$q\times q$\textit{\ submatrices ( with }$q=3$%
\textit{):}

\begin{equation*}
\mathit{\tilde{G}(y)}_{1}\mathit{=}\left[ 
\begin{array}{ccc}
2y_{1} & 0 & 3y_{3}^{2} \\ 
0 & 2y_{2} & 0 \\ 
2y_{1} & 2y_{2} & 0%
\end{array}%
\right] \mathit{,}\det \left( \tilde{G}(y)_{1}\right) \mathit{=-12y}_{1}%
\mathit{y}_{2}\mathit{y}_{3}^{2}
\end{equation*}

\begin{equation*}
\tilde{G}(y)_{2}=\left[ 
\begin{array}{ccc}
2y_{1} & 0 & 0 \\ 
0 & 2y_{2} & 3y_{4}^{2} \\ 
2y_{1} & 2y_{2} & 0%
\end{array}%
\right] ,\det \left( \tilde{G}(y)_{2}\right) =-12y_{1}y_{2}y_{4}^{2}
\end{equation*}

\begin{equation*}
\mathit{\tilde{G}(y)}_{3}\mathit{=}\left[ 
\begin{array}{ccc}
2y_{1} & 3y_{3}^{2} & 0 \\ 
0 & 0 & 3y_{4}^{2} \\ 
2y_{1} & 0 & 0%
\end{array}%
\right] \mathit{,}\det \left( \tilde{G}(y)_{3}\right) \mathit{=18y}_{1}%
\mathit{y}_{3}^{2}\mathit{y}_{4}^{2}
\end{equation*}

\begin{equation*}
\mathit{\tilde{G}(y)}_{4}\mathit{=}\left[ 
\begin{array}{ccc}
0 & 3y_{3}^{2} & 0 \\ 
2y_{2} & 0 & 3y_{4}^{2} \\ 
2y_{2} & 0 & 0%
\end{array}%
\right] \mathit{,}\det \left( \tilde{G}(y)_{4}\right) \mathit{=18y}_{2}%
\mathit{y}_{3}^{2}\mathit{y}_{4}^{2}
\end{equation*}

\textit{Hence }$\bar{\alpha}=4$\textit{\ but }$\det \left( \tilde{G}%
(y)_{3}\right) $\textit{\ and }$\det \left( \tilde{G}(y)_{4}\right) $\textit{%
\ are homogeneous polynomials of degree }$5>\bar{\alpha}.$

Thus, $\bar{\alpha}$\ is the smallest possible degree of an homogeneous
polynomial in the determinant of any non-singular $q\times q$\ submatrix of $%
G\left( y\right) $. Then $\bar{\alpha}=0$\ if and only if $y=0$ is not a
root of some such determinant and $\bar{\alpha}>0$\ otherwise. In other
words, $\bar{\alpha}=0$\ if and only if $G(0)$\ is of full row rank.

Select some matrix $\tilde{G}\left( y\right) _{l}$ for which $\bar{k}_{\det (%
\tilde{G}(y)_{l})}=\bar{\alpha}.$ Note that then $\left( \ref{homog}\right) $
\ implies that the limit:

\begin{equation}
\lim_{\lambda \rightarrow \infty }\lambda ^{\bar{\alpha}\ }\det \left( 
\tilde{G}_{l}(y/\lambda )\right)  \label{deg}
\end{equation}%
is a polynomial in $y$ on $R^{p}$ that is distinct from zero almost
everywhere.

For the matrix of polynomials $G\left( y\right) $ of rank $q$ and any
non-singular $q\times q$ matrix $S,$ for the polynomial matrix $SG\left(
y\right) $ there is some $\alpha =\left( \alpha _{1},...,\alpha _{p}\right) $
such that 
\begin{equation}
\lim_{\lambda \rightarrow \infty }diag(\lambda ^{\alpha _{i}})SG(y/\lambda )
\label{lim}
\end{equation}%
exists and is a finite non-zero polynomial matrix, $\bar{G}\left( y\right) .$
Indeed, define $\alpha _{i}=\underset{j}{\min \{}\bar{k}_{\left\{ SG\left(
y\right) \right\} _{ij}}\},$ where $\left\{ SG\left( y\right) \right\} _{ij}$
denotes the polynomial that is the $ij-th$ element of the matrix $SG\left(
y\right) .$ From $\left( \ref{homog}\right) $ existence of the limit matrix
follows.

\textbf{Lemma 4.1. }\textit{Suppose that there exists }$a=(\alpha
_{1},...,\alpha _{q})$\textit{\ with }$\alpha _{i}\geq 0$\textit{\ and a
non-singular }$q\times q$\textit{\ matrix }$S$\textit{\ such that\ the limit
matrix:}%
\begin{equation}
\bar{G}(y)=\lim_{\lambda \rightarrow \infty }diag(\lambda ^{\alpha
_{i}})SG(y/\lambda )  \label{g-bar}
\end{equation}%
\textit{is a finite non-zero matrix. Then for }$\bar{\alpha}$\textit{\ for
which (\ref{deg}) holds we get }$\sum\limits_{i=1}^{q}\alpha _{i}\leq \bar{%
\alpha}.$\textit{\ }$\bar{G}\left( y\right) $\textit{\ is non-singular if
and only if}

\begin{equation*}
\sum\limits_{i=1}^{q}\alpha _{i}=\bar{\alpha}.
\end{equation*}

\bigskip When $\bar{G}(y)$ exists for some $a$ and some matrix $S$, the
matrix $S$ can always be chosen such that $0\leq \alpha _{1}\leq ...\leq
\alpha _{q}\leq \bar{\alpha}.$

\textbf{Definition 4.1. }\textit{A }$q\times p$\textit{\ matrix of
polynomials }$G\left( y\right) $\textit{\ satisfies the "continuity of lower
degree ranks property" (CLDR) if for some non-singular }$q\times q$\textit{\
matrix }$S$\textit{\ and for some }$\alpha =(\alpha _{1},...,\alpha _{q})$%
\textit{\ such that }$\sum\limits_{i=1}^{q}\alpha _{i}=\bar{\alpha}$\textit{%
, }$0\leq \alpha _{1}\leq ...\leq \alpha _{q}\leq \bar{\alpha},$\textit{\ }$%
\left( \ref{g-bar}\right) $\textit{\ provides a rank }$q$\textit{\ matrix of
polynomials }$\bar{G}\left( y\right) .$

Essentially, the CLDR property holds if for some $S$ the transformed $%
SG\left( y\right) $ is such that the stabilizing rate $\bar{a}$ for the
determinant is shared between the rows of the matrix $SG\left( y\right) $
according to (\ref{g-bar}), and the limit matrix is non-singular.

The matrix $\bar{G}(y)$ depends upon the choice of the matrix $S$. Indeed in
Example 4.1 $\bar{\alpha}=4$ but it is clear that for $S=I$ Lemma 4.1 does
not hold. This is a consequence of the fact that there is a linear
dependence between the degree one polynomial terms in the rows of the
matrix. However setting

\begin{equation*}
S=\left[ 
\begin{array}{ccc}
1 & 0 & 0 \\ 
0 & 1 & 0 \\ 
1 & 1 & -1%
\end{array}%
\right]
\end{equation*}%
yields $SG\left( y/\lambda \right) $ as 
\begin{equation*}
\left[ 
\begin{array}{cccc}
2y_{1}/\lambda & 0 & 3y_{3}^{2}/\lambda ^{2} & 0 \\ 
0 & 2y_{2}/\lambda & 0 & 3y_{4}^{2}/\lambda ^{2} \\ 
0 & 0 & 3y_{3}^{2}/\lambda ^{2} & 3y_{4}^{2}/\lambda ^{2}%
\end{array}%
\right]
\end{equation*}%
and CLDR holds with this $S$ and $\alpha =(1,1,2)$.

The next example demonstrates that the CLDR property may not hold for some $%
G\left( y\right) $ even with $q=p.$

\textbf{Example 4.2. }\textit{Consider}

\begin{equation*}
G(y)=\left[ 
\begin{array}{cc}
y_{1} & 0 \\ 
(c+y_{2})^{2} & y_{1}(c+y_{2})%
\end{array}%
\right]
\end{equation*}%
\textit{with }$c\neq 0.$\textit{\ Then }$\bar{\alpha}=2$\textit{. Consider
an arbitrary }$2\times 2$\textit{\ matrix }$S=(s_{ij})_{1\leq i,j\leq 2}.$

\textit{Three possibilities could arise for }$\bar{\alpha}=2$\textit{\ if
CLDR were to hold so that }$\bar{a}=\alpha _{1}+\alpha _{2}.$

\textit{First, }$\alpha _{1}=\alpha _{2}=1,$\textit{\ then}%
\begin{equation*}
\lim_{\lambda \rightarrow \infty }\left[ 
\begin{array}{cc}
\lambda & 0 \\ 
0 & \lambda%
\end{array}%
\right] SG(y/\lambda )
\end{equation*}%
\textit{does not exist, except if }$s_{12}=s_{22}=0$\textit{, which is
precluded for non-singular matrix }$S$\textit{.}

\textit{Second, }$\alpha _{1}=2,$\textit{\ }$\alpha _{2}=0,$\textit{\ then}%
\begin{equation*}
\lim_{\lambda \rightarrow \infty }\left[ 
\begin{array}{cc}
\lambda ^{2} & 0 \\ 
0 & 1%
\end{array}%
\right] SG(y/\lambda )
\end{equation*}%
\textit{does not exist for any non-zero matrix }$S$\textit{.}

\textit{Third, }$\alpha _{1}=0,$\textit{\ }$\alpha _{2}=2,$\textit{\ then}%
\begin{equation*}
\lim_{\lambda \rightarrow \infty }\left[ 
\begin{array}{cc}
1 & 0 \\ 
0 & \lambda ^{2}%
\end{array}%
\right] SG(y/\lambda )
\end{equation*}%
\textit{does not exist, except if }$s_{21}=s_{22}=0$\textit{, which is
precluded for non-singular matrix }$S$\textit{.}

Now that we see that some matrices of polynomials satisfy the CLDR property
and some do not, we further characterize the difference between the two
possibilities.

\textbf{Lemma 4.2. }\textit{Given a matrix }$G\left( y\right) $\textit{\
with the corresponding }$\bar{a},$\textit{\ for any non-singular matrix }$S$%
\textit{\ and }$a^{\prime }=\left( \alpha _{1}^{\prime },...,\alpha
_{q}^{\prime }\right) $\textit{\ with }$0\leq \alpha _{1}^{\prime }\leq
...\leq \alpha _{q}^{\prime }\leq \bar{\alpha}$\textit{\ and }$%
\sum\limits_{i=1}^{q}\alpha _{i}^{\prime }=\bar{\alpha}$\textit{\ either (i)
CLDR property holds with this }$S$\textit{\ and }$a^{\prime },$\textit{\ or
(ii) no finite limit exists for}%
\begin{equation*}
\left[ diag(\lambda ^{\alpha _{i}^{\prime }})SG(y/\lambda )\right] ,
\end{equation*}%
\textit{or (iii) if a finite limit does exist }%
\begin{equation}
rank\lim_{\lambda \rightarrow \infty }\left[ diag(\lambda ^{\alpha
_{i}^{\prime }})SG(y/\lambda )\right] <q.  \label{defrank}
\end{equation}

Thus if $S$ and $a$ are such that a finite limit $\left( \ref{g-bar}\right) $
exists then either the CLDR property holds for such $S,a$ or the limit
matrix $\bar{G}\left( y\right) $ has a deficient rank. If the limit matrix $%
\bar{G}\left( y\right) $ has a deficient rank for some $S,a,$ it has a
deficient rank for any other \thinspace $S^{\prime },a^{\prime }.$ We can
thus say that $G\left( y\right) $ is either CLDR or deficient rank. To
determine whether there exist some $S$ and $a$ for which CLDR property holds
we provide a recursive construction of $S$ and $a$ that either gives the
CLDR property or results in a deficient rank.

\textbf{Lemma 4.3. }\textit{Given a }$q\times p$\textit{\ matrix }$G\left(
y\right) $\textit{\ of polynomials, there is a recursive construction that
provides the pair }$S$\textit{\ and }$a,$\textit{\ such that either CLDR
property is satisfied for this pair or the deficient rank property holds.}

The construction in the proof implies that we can write:

\begin{equation}
SG(y)=\bar{G}(y)+\bar{R}(y)  \label{decomp}
\end{equation}%
where for $i=1,...,q,$ the row $i$\ of $\bar{R}(y)$\ contains no homogeneous
polynomial of order smaller or equal to $\alpha _{i}.$

\subsection{Vectors of polynomial functions, Jacobian matrices and the Wald
statistic}

Consider the $q\times 1$ vector of polynomial functions, $g\left( y\right) $
with $g\left( 0\right) =0$ and the Jacobian matrix of polynomials, $G(y)=%
\frac{\partial g}{\partial y^{\prime }}(y).$

Consider a non-singular $S$ that satisfies $\left( \ref{decomp}\right) $ for 
$G\left( y\right) $ (and $\sum\limits_{i=1}^{q}\alpha _{i}\leq \bar{\alpha}%
). $

Then 
\begin{equation}
Sg(y)=\bar{g}(y)+\bar{r}(y)  \label{decomp1}
\end{equation}%
where for every $i=1,...,q$

\begin{eqnarray*}
\bar{g}_{i}(y) &=&\int_{0}^{y}\bar{G}(x)_{i}dx; \\
\bar{r}_{i}(y) &=&\int_{0}^{y}\bar{R}(x)_{i}dx,
\end{eqnarray*}%
where the integration of the gradient along any continuous curve from $0$ to 
$y$ provides each component of $g,r.$

\bigskip Each $\bar{g}_{i}(y)$\ of\ $\bar{g}(y)$\ is a homogeneous
polynomial of order $(\alpha _{i}+1)$ and, by Euler formula:%
\begin{equation}
\bar{g}(y)=\Lambda \bar{G}(y)y  \label{Euler}
\end{equation}%
with:

\begin{equation*}
\Lambda =diag\left( \frac{1}{\alpha _{i}+1}\right) .
\end{equation*}%
Each element $\bar{r}_{i}(y)$\ of $\bar{r}(y)$\ contains no homogeneous
polynomial of order smaller or equal to $(\alpha _{i}+1).$

In particular, when $\lambda $\ goes to infinity:

\begin{eqnarray}
diag(\lambda ^{\alpha _{i}})SG(y/\lambda ) &=&\bar{G}(y)+O(1/\lambda );
\label{decomp2} \\
diag(\lambda ^{\alpha _{i}})S\lambda g(y/\lambda ) &=&\bar{g}(y)+O(1/\lambda
).  \notag
\end{eqnarray}

\bigskip Define now for some positive definite matrix $\Omega $ a quadratic
form 
\begin{equation}
W(y,g,\lambda ,\Omega )=\lambda ^{2}g^{\prime }(y/\lambda )[G(y/\lambda
)\Omega G^{\prime }(y/\lambda )]^{-1}g(y/\lambda ).  \label{quad}
\end{equation}

Note that $W\left( y,g,\lambda ,\Omega \right) =W\left( y,Mg,\lambda ,\Omega
\right) $ for any non-singular matrix $M;$ we can choose $M=S(\lambda
)=diag(\lambda ^{\alpha _{i}})S.$ This provides

\begin{equation}
W(y,g,\lambda ,\Omega )=g^{\prime }(y/\lambda )\lambda S^{\prime }(\lambda
)[S(\lambda )G(y/\lambda )\Omega G^{\prime }(y/\lambda )S^{\prime }(\lambda
)]^{-1}S(\lambda )\lambda g(y/\lambda ).  \label{asymp}
\end{equation}

Suppose that $\Omega =\Omega \left( \lambda \right) $ with the property that
as $\lambda \rightarrow \infty $ the matrix $\Omega =\Omega ^{0}+o(1),$ with 
$\Omega ^{0}$ a non-singular matrix. Then we can write as $\lambda
\rightarrow \infty $

\begin{eqnarray*}
&&W(y,g,\lambda ,\Omega ) \\
&=&[\bar{g}(y)+O(1/\lambda )]^{\prime }\left\{ \left[ \bar{G}(y)+O(1/\lambda
)\right] [\Omega ^{0}+o(1)]\left[ \bar{G}(y)+O(1/\lambda )\right] ^{\prime
}\right\} ^{-1}[\bar{g}(y)+O(1/\lambda )].
\end{eqnarray*}

If CLDR property holds for $G,$ $\bar{G}\left( y\right) $ is full rank and
then

\begin{eqnarray}
\lim_{\lambda \rightarrow \infty }W(y,g,\lambda ,\Omega ) &=&[\bar{g}%
(y)]^{\prime }\left\{ \left[ \bar{G}(y)\right] \Omega ^{0}\left[ \bar{G}(y)%
\right] ^{\prime }\right\} ^{-1}[\bar{g}(y)]  \label{limit} \\
&=&W_{\infty }(y,g,\Omega ^{0}).  \notag
\end{eqnarray}

\bigskip Next, we demonstrate tht if CLDR property does not hold $%
W(y,g,\lambda ,\Omega )$ diverges to infinity as $\lambda \rightarrow \infty
.$

Suppose that CLDR property does not hold, then find $a$ for which $\left( %
\ref{lim}\right) $ provides a finite matrix, by lack of CLDR in that case $%
\Sigma \alpha _{i}<\bar{a}.$

Then recall that $[G(y/\lambda )\Omega G^{\prime }(y/\lambda )]^{-1}$ can be
represented as the ratio of the adjoint matrix, denoted $[G(y/\lambda
)\Omega G^{\prime }(y/\lambda )]^{\ast },$ to the determinant, $\det
[G(y/\lambda )\Omega G^{\prime }(y/\lambda )].$ Write (\ref{asymp}) as%
\begin{equation*}
\frac{g^{\prime }(y/\lambda )\lambda S^{\prime }(\lambda )[S(\lambda
)G(y/\lambda )\Omega G^{\prime }(y/\lambda )S^{\prime }(\lambda )]^{\ast
}S(\lambda )\lambda g(y/\lambda )}{\det [S(\lambda )G(y/\lambda )\Omega
G^{\prime }(y/\lambda )S^{\prime }(\lambda )]};
\end{equation*}%
this is 
\begin{equation*}
\frac{\lbrack \bar{g}(y)+O(1/\lambda )]^{\prime }\left\{ \left[ \bar{G}%
(y)+O(1/\lambda )\right] [\Omega ^{0}+o(1)]\left[ \bar{G}(y)+O(1/\lambda )%
\right] ^{\prime }\right\} ^{\ast }[\bar{g}(y)+O(1/\lambda )]}{\det
[S(\lambda )G(y/\lambda )\Omega G^{\prime }(y/\lambda )S^{\prime }(\lambda )]%
}.
\end{equation*}%
The numerator has a finite limit.

In the denominator we have 
\begin{eqnarray*}
\det [S(\lambda )G(y/\lambda )\Omega G^{\prime }(y/\lambda )S^{\prime
}(\lambda )] &=&\lambda ^{2\Sigma \alpha _{i}}\det [SG(y/\lambda )\Omega
G^{\prime }(y/\lambda )S^{\prime }] \\
&=&\lambda ^{2[\Sigma \alpha _{i}-\bar{\alpha}]}\lambda ^{2\bar{\alpha}}\det
[SG(y/\lambda )\Omega G^{\prime }(y/\lambda )S^{\prime }].
\end{eqnarray*}%
Thus as $\lambda \rightarrow \infty ,$ when the CLDR property is not
fulfilled, while $\lambda ^{2\bar{\alpha}}\det [SG(y/\lambda )\Omega
G^{\prime }(y/\lambda )S^{\prime }]$ has a finite limit for every $\Omega $, 
$\lambda ^{2[\Sigma \alpha _{i}-\bar{\alpha}]}$ converges to zero and 
\begin{equation}
W(y,g,\lambda ,\Omega )\underset{\lambda \rightarrow \infty }{\rightarrow }%
\infty .  \label{blowup}
\end{equation}

Thus CLDR property plays a very important role in the existence of a limit
for the Wald statistic.

\subsection{The limit distribution of the Wald statistic}

Define $y=A(\theta -\bar{\theta})$ for some non-degenerate matrix $A;$ with
this substitution under the assumption $g\left( \bar{\theta}\right) =0$ the
polynomial function $g\left( \theta \right) \,\ $becomes $g\left( A^{-1}y+%
\bar{\theta}\right) =g_{\bar{\theta}}(y),$ a polynomial function with $g_{%
\bar{\theta}}(0)=0.$ The Jacobian polynomial matrix $G$ gets multiplied by
the nonsingular matrix $A^{-1}$ to provide the new Jacobian $G_{\bar{\theta}%
}\left( y\right) $ with respect to $y$ and $G\left( \theta \right) =G_{\bar{%
\theta}}(y)A.$ Note the role that the nonsingular matrix $A$ plays: it does
not change the order of polynomial function $g;$ if CLDR property holds for $%
G_{\bar{\theta}}$ defined with some nonsingular $A,$ it holds for any other
nonsingular $A.$ In the notation for the function and the Jacobian we do not
emphasize then the role of $A.$

The Wald test statistic in $\left( \ref{waldgen}\right) $ for $\theta =\hat{%
\theta}_{T},\lambda =\lambda _{T},\Omega =\hat{V}_{T}$ and with $%
y_{T}=A\left( \hat{\theta}_{T}-\bar{\theta}\right) $ can be written as

\begin{equation*}
\lambda _{T}^{2}g_{\bar{\theta}}^{\prime }(y_{T})[G_{\bar{\theta}}(y_{T})A%
\hat{V}_{T}A^{\prime }G_{\bar{\theta}}^{\prime }(y_{T})]^{-1}g_{\bar{\theta}%
}(y_{T}).
\end{equation*}%
Consider $g_{\bar{\theta},0}(y)$ for $A_{0}=V^{-\frac{1}{2}},$ and $g_{\bar{%
\theta},1}(y)$ for some nonsingular $A;$ then the distribution of $g_{\bar{%
\theta},1}(V^{\frac{1}{2}}AZ)$ is the same as $g_{\bar{\theta},0}(Z).$ The
distribution of $G_{\bar{\theta},1}(V^{\frac{1}{2}}AZ)$ is the same as $G_{%
\bar{\theta},0}(Z)V^{-\frac{1}{2}}A^{-1},$ since we consider non-singular
reparametrizations of continuous functions. Thus the limit distribution does
not depend on $A.$ Below we write $g_{\bar{\theta}},$ $G_{\bar{\theta}}$ for 
$g_{\bar{\theta},0}$ and $G_{\bar{\theta},0}.$

\textbf{Theorem 4.1. }\textit{Under the Assumptions 2.1a, 2.2, 2.3 and 2.4
if (a) at }$\bar{\theta}$\textit{\ the CLDR property holds for }$G_{\bar{%
\theta}}(y)$\textit{\ (for any non-singular }$A)$\textit{\ then the limit
distribution of }$W_{T}$\textit{\ as }$T\rightarrow \infty $\textit{\ is
given by the distribution of }%
\begin{equation}
\lbrack \bar{g}_{\bar{\theta}}(Z)]^{\prime }\left\{ \left[ \bar{G}_{\bar{%
\theta}}(Z)\right] [\bar{G}_{\bar{\theta}}(Z)]^{\prime }\right\} ^{-1}[\bar{g%
}_{\bar{\theta}}(Z)];  \label{genwaldlim}
\end{equation}%
\textit{if (b) at }$\bar{\theta}$\textit{\ the deficient rank property
holds, then }$W_{T}$\textit{\ diverges to infinity as }$T\rightarrow \infty
. $

\textbf{Corollary 4.1. }\textit{When the CLDR property holds the limit
distribution can be represented as the distribution of}%
\begin{equation}
Z^{\prime }[\bar{G}_{\bar{\theta}}(Z)]^{\prime }\Lambda _{\bar{\theta}%
}\left\{ \left[ \bar{G}_{\bar{\theta}}(Z)\right] [\bar{G}_{\bar{\theta}%
}(Z)]^{\prime }\right\} ^{-1}\Lambda _{\bar{\theta}}[\bar{G}_{\bar{\theta}%
}(Z)]Z.  \label{eulerwaldlim}
\end{equation}

This follows from the Euler formula $\left( \ref{Euler}\right) .$

The Example below illustrates the applicability of parts (a) and (b) of
Theorem 4.1.

\textbf{Example 4.3. }\textit{Recall Example 3.4 with }$g(\theta )=\left[ 
\begin{array}{c}
\theta _{1}^{2} \\ 
\theta _{1}\theta _{2}^{2}%
\end{array}%
\right] .$\textit{\ Then, the set of possible values of }$\bar{\theta}$%
\textit{\ under the null is the line }$(0,\bar{\theta}_{2})^{\prime },$%
\textit{\ }$\bar{\theta}_{2}\in 
\mathbb{R}
,$\textit{\ and for }$A=I$

\begin{equation*}
\mathit{G(y)=}\left[ 
\begin{array}{cc}
2y_{1} & 0 \\ 
(\bar{\theta}_{2}+y_{2})^{2} & 2y_{1}(\bar{\theta}_{2}+y_{2})%
\end{array}%
\right] \mathit{.}
\end{equation*}

\textit{When }$\bar{\theta}_{2}\neq 0$\textit{\ }$\bar{G}=\left[ 
\begin{array}{cc}
2y_{1} & 0 \\ 
\bar{\theta}_{2}^{2} & 0%
\end{array}%
\right] ;$\textit{\ the CLDR property does not hold, (b) of the Theorem
applies. By contrast, if }$\bar{\theta}_{2}=0$\textit{, we have }$\bar{\alpha%
}=3,$\textit{\ and the sharing rule }$\alpha =(1,2)$\textit{\ for CLDR
immediately follows and (a) applies.}

When there is only one restriction there is only one $\alpha _{1}=\alpha _{%
\bar{\theta}}.$ Here CLDR always holds and thus under the assumptions of
Theorem 4.1 the convergence of $W_{T}$ to%
\begin{equation*}
\frac{1}{\left( 1+\alpha _{\bar{\theta}}\right) ^{2}}\frac{\left( Z^{\prime }%
\bar{G}_{\bar{\theta}}\left( Z\right) ^{\prime }\right) ^{2}}{\bar{G}_{\bar{%
\theta}}\left( Z\right) \bar{G}_{\bar{\theta}}\left( Z\right) ^{\prime }}%
\equiv \frac{\left\Vert \bar{g}\left( Z\right) \right\Vert ^{2}}{\left\Vert 
\bar{G}\left( Z\right) \right\Vert ^{2}}
\end{equation*}%
always obtains.

In the case of multiple restrictions violation of the CLDR property is
possible; in such a case the statistic may diverge under the null. One could
consider replacing the restrictions by a set of equivalent restrictions that
preclude violation of CLDR property. This is always possible.

Indeed, for any vector $g$ of $q$ restrictions $g\left( \theta \right) =0,$
the restrictions are equivalent to a single restriction 
\begin{equation*}
\left\Vert g(\theta )\right\Vert ^{2}=\sum\limits_{i=1}^{q}g_{i}^{2}(\theta
)=0.
\end{equation*}%
Since the CLDR property is not an issue with one constraint a possible
strategy is to replace the $q$ restriction by the single restriction and
consider the corresponding test statistic.

Of course, this simplification may have an important cost in terms of power
since itdoes not take into account the fact that for an estimator $\hat{%
\theta}_{T}$ the components may be highly correlated. Then, the naive norm $%
\left\Vert g(\hat{\theta}_{T})\right\Vert $\ of the vector $g(\hat{\theta}%
_{T})$\ may not be the efficient way to assess its distance from zero; some
weighting may be advantageous.

\section{Bounds on the statistic and bounds on critical values}

Sometimes the asymptotic distribution of the Wald statistic under the null,
even when non-standard, can be uniquely determined; this is the case in
Example 3.3. But typically under conditions of Theorem 4.1 with possible
singularity the asymptotic distribution under the null is not uniquely
determined.\ Because the asymptotic distribution of the Wald statistic may
be discontinuous in the true values it is useful to establish uniform bounds
on the asymptotic distribution of the statistic, or on the critical values
for the test.

Denote by $\alpha $ the smallest $\alpha _{i}$ (usually $\alpha _{1})$ in
Definition 4.1. Below we show that $\frac{1}{\left( 1+\alpha \right) ^{2}}%
\left\Vert Z\right\Vert ^{2}$ (distributed $\frac{1}{\left( 1+\alpha \right)
^{2}}\chi _{p\text{ }}^{2}$ under normality) provides a uniform upper bound
on the asymptotic null-distribution always under conditions of the Theorem
4.1(a), i.e. when CLDR property holds. If $\alpha =0,$ then there may be no
singularity in which case with normality the usual asymptotic $\chi _{q}^{2}$
distribution holds; in general the overall uniform bound with or without
singularity under Theorem 4.1.$\left( a\right) $ is $\chi _{p}^{2}$.

It is possible to improve on the $\chi _{p}^{2}$ bound when $\alpha \geq 1$;
sometimes the form of the restrictions may provide $\alpha \geq 1$. When
this is not the case, it may be possible to establish that $\alpha \geq 1$
by the adaptive strategy proposed in the next section that would eliminate
the possibility that $\alpha =0.$

However, for testing it may be sufficient to bound the distribution in the
tail rather than everywhere, and so uniform bounds on critical values are
also of interest. Gouri\'{e}roux and Jasiak (2013) discuss a bound on
critical values for a test of a determinant.

Here in Theorem 5.1 we first establish general bounds on an asymptotic
distribution derived for a particular vector of true parameter values, when
there may be a singularity at that value. We separately examine the case of
one restriction. We also examine a relation between critical values at
different $\alpha .$ In the case of one restriction it is possible to
provide the number of variables for which generally the standard critical
values deliver a conservative test.

\subsection{A general uniform upper bound}

Start with the representation of the asymptotic distribution from (\ref%
{eulerwaldlim}): 
\begin{equation*}
W\left( Z\right) =Z^{\prime }\bar{G}(Z)^{\prime }\Lambda \left[ \bar{G}(Z)%
\bar{G}(Z)^{\prime }\right] ^{-1}\Lambda \bar{G}(Z)Z.
\end{equation*}

This distribution depends on the singularity properties that are exhibited
at the true value $\bar{\theta},V.$

As the theorem below states a bound that depends only on $\alpha $ is
possible in all cases when CLDR holds.

\bigskip \textbf{Theorem 5.1. }\textit{Under the conditions of Theorem
4.1(a), the asymptotic distribution of the Wald statistic under the null
(that depends on the singularity properties at }$\bar{\theta})$\textit{\ is
bounded from above by the distribution of }$\frac{1}{\left( 1+\alpha \right)
^{2}}\left\Vert Z\right\Vert ^{2};$\textit{\ under the normality Assumption
2.3 this bound is }$\frac{1}{\left( 1+\alpha \right) ^{2}}\chi _{p}^{2}.$

Thus under conditions of the Theorem 5.1 there is always a general upper
bound on the distribution of the Wald statistic under the null given by $%
\chi _{p}^{2}.$

\textbf{Remark 5.1. }\textit{When }$\Lambda =I$\textit{\ implying that }$%
\bar{G}(Z)$\textit{\ does not depend on }$Z$\textit{\ and is a }$q\times p$%
\textit{\ rank }$q$\textit{\ matrix of constants the projection is onto a }$%
q-$\textit{dimensional subspace, the limit distribution is standard and is
given under normality by }$\chi _{q}^{2}.$

\bigskip In the case of one restriction under normality the upper bound is
either the usual $\chi _{1}^{2},$ if $\alpha =0,$ or else for some $\alpha
>0 $ the bound is $\frac{1}{(1+\alpha )^{2}}\chi _{p}^{2}.$ If all that is
known is that $\alpha >0,$ then the bound $\frac{1}{4}\chi _{p}^{2}$ applies
for any such $\alpha .$

In Example 3.3 the limiting distribution is $\frac{1}{4}\chi _{p}^{2};$ thus
this bound can be attained.

In the special case $p=q$ and under CLDR $\bar{G}$ is invertible a.e. and $%
Z^{\prime }\bar{G}(Z)^{\prime }\Lambda \left[ \bar{G}(Z)\bar{G}(Z)^{\prime }%
\right] ^{-1}\Lambda \bar{G}(Z)Z=Z^{\prime }\bar{G}^{\prime }\Lambda \bar{G}%
^{\prime -1}\bar{G}^{-1}\Lambda \bar{G}Z.$ Then 
\begin{equation*}
W(Z)\leq \left\Vert \Lambda \right\Vert ^{2}\left\Vert Z\right\Vert ^{2},
\end{equation*}%
since the norm of a similar matrix is the same as for $\Lambda .$ Under
normality the bound is $\frac{1}{\left( 1+\alpha \right) ^{2}}\chi _{q}^{2}.$
In this case the asymptotic distribution is bounded from above by the usual
distribution and under normality the distribution $\chi _{q}^{2}$ provides a
conservative test.

\subsection{Bounds on critical values for purely singular cases $\protect%
\alpha \geq 1$ under normality}

A conservative test for a given level may be given by the standard critical
values, even in the non-standard cases considered here since then dominance
by the standard distribution is required only in the tail and not
everywhere. The following Lemma demonstrates that when the distribution is
purely singular ($\alpha \geq 1)$ there is always a level, $\gamma _{0},$
such that using the standard critical values for any $\gamma \leq \gamma
_{0} $ provides a conservative asymptotic test. Indeed, there exists $\gamma
_{0}$ such that $\Pr (\frac{1}{\left( 1+\alpha \right) ^{2}}\chi
_{p}^{2}>\chi _{q}^{2}(\gamma _{0}))<\gamma _{0},$ where $\chi
_{q}^{2}\left( \gamma _{0}\right) $ denotes the critical value.

The following lemma establishes this tail dominance.

\textbf{Lemma 5.1. }\textit{Consider two random variables }$T\sim \chi
_{p_{1}}^{2}/\alpha _{1}$\textit{\ and }$S\sim \chi _{p_{2}}^{2}/\alpha
_{2}, $\textit{\ where }$p_{2}>p_{1},\alpha _{2}>\alpha _{1}>0.$\textit{Then
there exists }$y_{0}$\textit{\ such that for }$y>y_{0}$\textit{\ we have}

\begin{equation*}
p.d.f._{S}(y)<p.d.f._{T}(y).
\end{equation*}

This makes it possible to rely only on $p$ and $q$ in indicating when
standard critical values provide a conservative test.

When there may be a singularity with $\alpha \geq 1$ the critical value
coming from the standard test will at some level result in a conservative
Wald test; the question is whether this holds for conventional test levels.
Abstracting from the specific form of restrictions the answer depends on $%
\alpha ,$ $p$ and $\dot{q};$ the \ higher the $\alpha $ and the closer
together $p$ and $q,$ the easier to obtain conservative asymptotic tests at
conventional levels. Since $p$ and $q$ are given by the restrictions, all
that is required is to establish $\alpha .$

Comparing the values of $p.d.f._{\chi _{q}^{2}}(y_{.05})$ with $p.d.f._{\chi
_{p}^{2}/\left( 1+a\right) ^{2}}(y_{.05})$ where $y_{.05}$ is the critical
value for $\chi _{q}^{2}$ at $.05$ level we determine for which $\max p$ we
get a smaller value for the second $p.d.f.;$ because of monotonicity in the
tail this indicates smaller probability and a conservative test.

For one restriction the standard test based on $\chi _{1}^{2}$ critical
value is conservative at $.05$ level for $p\leq 6$ but may not be not for $%
p=7.$ At $.01$ level this test is conservative for $p\leq 10,$ but may not
be for $p=11.$ To show this only a computation of the critical values for $%
\chi _{1}^{2}$ and for overall bound $\frac{1}{4}\chi _{p}^{2}$ is required.

When CLDR holds for $q=2,$ if $\alpha =1$ at $.05$ level we get $\max p=11,$
for $q=3$ and $\alpha =1$ we get $\max p=17.$

These computations show that in many situations the standard test is
conservative.

\section{An adaptive strategy for determining the asymptotic distribution
and the bounds}

From $\left( \ref{eulerwaldlim}\right) $ it follows that the asymptotic
distribution requires the knowledge of $\bar{G}_{\bar{\theta}}$ and $\Lambda
_{\bar{\theta}}.$ For bounds determining the lowest value on the diagonal of 
$\Lambda _{\bar{\theta}}$ is sufficient. The construction in proof of Lemma
4.3 makes it clear that the main issue for finding the elements of $\Lambda
_{\bar{\theta}}$ is deciding on the lowest order of the homogeneous
polynomial that enters non-trivially into a polynomial (that represents a
matrix entry or a determinant of a polynomial matrix), to find $\bar{G}_{%
\bar{\theta}}$ homogeneous polynomials (their coefficients) of the
corresponding lowest orders would have to be consistently estimated.

\subsection{Adaptive estimation of polynomial functions and orders}

For $\bar{\theta}$ consider $P_{\bar{\theta}}(\theta )$ that is a polynomial
function of order $m_{P}$ in components of a $p\times 1$\ vector $\bar{\theta%
}$ with the representation in terms of components of $\theta -\bar{\theta}$
given by%
\begin{eqnarray}
P_{\bar{\theta}}(\theta ) &=&\bar{P}_{0}\left( \bar{\theta}\right)
+\tsum\limits_{k=1}^{m_{P}}\bar{P}_{k}(\theta -\bar{\theta})
\label{genericpolyn} \\
&=&\bar{P}\left( 0,...,0,\bar{\theta}\right) +\tsum\limits_{k=1}^{m_{P}}%
\left[ \tsum\limits_{i_{1}+...+i_{p}=k}\bar{P}\left( i_{1},...,i_{p},\bar{%
\theta}\right) \left( \theta _{1}-\bar{\theta}_{1}\right) ^{i_{1}}...\left(
\theta _{p}-\bar{\theta}_{p}\right) ^{i_{p}}\right] ,  \notag
\end{eqnarray}%
where the corresponding coefficients $\bar{P}\left( i_{1},...,i_{p},\bar{%
\theta}\right) $ are values of a polynomial in components of $\bar{\theta}$
and the constant term $\bar{P}_{0}\left( \bar{\theta}\right) $ can be
represented as a coefficient, $\bar{P}\left( 0,...,0,\bar{\theta}\right) ;$ $%
P_{\bar{\theta}}\left( \bar{\theta}\right) =\bar{P}_{0}\left( \bar{\theta}%
\right) .$

Consider a linear substitution with a nonsingular matrix $A:$ $y=A\left(
\theta -\bar{\theta}\right) ,$ with it 
\begin{equation*}
\left( \theta _{1}-\bar{\theta}_{1}\right) ^{i_{1}}...\left( \theta _{p}-%
\bar{\theta}_{p}\right) ^{i_{p}}=\tsum\limits_{\substack{ i_{1}^{\prime
}+...+i_{p}^{\prime }=k \\ 1\leq i_{v}^{\prime }\leq i_{v}}}\bar{A}%
(i_{1}^{\prime },...,i_{p}^{\prime },A)y_{1}^{i_{1}^{\prime
}}...y_{p}^{i_{p}^{\prime }},
\end{equation*}%
with some coefficients $\bar{A}(i_{1}^{\prime },...,i_{p}^{\prime },A)$ that
are polynomials in the matrix elements of the matrix $A^{-1}.$ Then the
polynomial $P_{\bar{\theta}}(\theta )$ becomes%
\begin{eqnarray}
P_{\bar{\theta}}(y) &=&\bar{P}\left( 0,...,0,\bar{\theta}\right)
+\tsum\limits_{k=1}^{m_{P}}\left[ \tsum\limits_{i_{1}+...+i_{p}=k}\bar{P}%
\left( i_{1},...,i_{p},\bar{\theta}\right) \tsum\limits_{\substack{ %
i_{1}^{\prime }+...+i_{p}^{\prime }=k \\ 1\leq i_{v}^{\prime }\leq i_{v}}}%
\bar{A}_{i_{1},...,i_{p}}(i_{1}^{\prime },...,i_{p}^{\prime
},A)y_{1}^{i_{1}^{\prime }}...y_{p}^{i_{p}^{\prime }}\right]  \\
&=&\bar{P}\left( 0,...,0,\bar{\theta}\right)
+\tsum\limits_{k=1}^{m_{P}}P_{k}\left( y\right)   \label{homog1} \\
&=&\bar{P}\left( 0,...,0,\bar{\theta}\right) +\tsum\limits_{k=1}^{m_{P}}%
\left[ \tsum\limits_{i_{1}^{\prime }+...+i_{p}^{\prime }=k}\left(
\tsum\limits_{\substack{ i_{1}+...+i_{p}=k \\ i_{v}^{\prime }\leq i_{v}}}%
\bar{A}_{i_{1},...,i_{p}}(i_{1}^{\prime },...,i_{p}^{\prime },A)\bar{P}%
\left( i_{1},...,i_{p},\bar{\theta}\right) \right) y_{1}^{i_{1}^{\prime
}}...y_{p}^{i_{p}^{\prime }}\right]   \notag \\
&=&\bar{P}\left( 0,...,0,\bar{\theta}\right) +\tsum\limits_{k=1}^{m_{P}}%
\left[ \tsum\limits_{i_{1}+...+i_{p}=k}\bar{P}\left( i_{1},...,i_{p},\bar{%
\theta},A\right) y_{1}^{i_{1}}...y_{p}^{i_{p}}\right] ,  \notag
\end{eqnarray}%
where the coefficients $\bar{P}\left( i_{1},...,i_{p},\bar{\theta},A\right)
=\tsum\limits_{\substack{ i_{1}^{\prime }+...+i_{p}^{\prime }=k \\ i_{v}\leq
i_{v}^{\prime }}}\bar{A}_{i_{1}^{\prime },...,i_{p}^{\prime
}}(i_{1},...,i_{p},A)\bar{P}\left( i_{1}^{\prime },...,i_{p}^{\prime },\bar{%
\theta}\right) .$

For estimator $\hat{\theta}_{T}$ of $\bar{\theta}$ define estimators of the
coefficients $\bar{P}\left( i_{1},...,i_{k},\bar{\theta}\right) $ by 
\begin{equation}
\hat{P}\left( i_{1},...,i_{p},\bar{\theta}\right) =\left\{ 
\begin{array}{cc}
\bar{P}\left( i_{1},...,i_{p},\hat{\theta}_{T}\right) & \text{if }\left\vert 
\bar{P}\left( i_{1},...,i_{p},\hat{\theta}_{T}\right) \right\vert \geq \frac{%
c}{\lambda _{T}^{\delta }}; \\ 
0 & \text{if }\left\vert \bar{P}\left( i_{1},...,i_{p},\hat{\theta}%
_{T}\right) \right\vert <\frac{c}{\lambda _{T}^{\delta }}%
\end{array}%
\right.  \label{estimator}
\end{equation}%
for $0<\delta <1$ and some $c>0.$

If $A=I,$ no further estimation is required.

For $A=V^{-\frac{1}{2}}$ and an estimator $\hat{V}_{T}$ of $V$ estimate $%
\bar{A}_{i_{1}^{\prime },...,i_{p}^{\prime }}(i_{1},...,i_{p},V^{-\frac{1}{2}%
})$ by $\hat{A}_{i_{1}^{\prime },...,i_{p}^{\prime }}(i_{1},...,i_{p},V^{-%
\frac{1}{2}})=\bar{A}_{i_{1}^{\prime },...,i_{p}^{\prime }}(i_{1},...,i_{p},%
\hat{V}_{T}^{-\frac{1}{2}}).$

Combining according to $\left( \ref{homog1}\right) $ we can obtain the
estimator $\hat{P}\left( i_{1},...,i_{p},\bar{\theta},A\right) $ of $\bar{P}%
\left( i_{1},...,i_{p},\bar{\theta},A\right) .$

Define (as in $\left( \ref{k-bar}\right) )$

\begin{equation}
k_{P}=\underset{0\leq k\leq m_{P}}{\min }\{k:\bar{P}\left( i_{1},...,i_{p},%
\bar{\theta},A\right) \neq 0\text{ for }i_{1}+...+i_{p}=k\},  \label{k-bar1}
\end{equation}%
and correspondingly $\hat{k}_{P}$ for the polynomial $P$ with estimated
coefficients $\hat{P}\left( i_{1},...,i_{p},\bar{\theta},A\right) .$ Note
that $k_{p}$ does not depend on $A.$

$\ $

\textbf{Lemma 6.1.} \textit{For }$\hat{\theta}_{T}$\textit{\ and }$\hat{V}%
_{T}$\textit{\ that satisfy Assumptions 2.3 and 2.4 }%
\begin{equation*}
\hat{P}\left( i_{1},...,i_{p},\bar{\theta},V\right) -\bar{P}\left(
i_{1},...,i_{p},\bar{\theta},V\right) \rightarrow _{p}0,
\end{equation*}%
\textit{moreover if }$\bar{P}\left( i_{1},...,i_{p},\bar{\theta}\right) =0,$%
\textit{\ }%
\begin{equation*}
\Pr \left( \hat{P}\left( i_{1},...,i_{p},\bar{\theta}\right) =0\right)
\rightarrow 1
\end{equation*}%
\textit{and} 
\begin{equation*}
\Pr \left( \hat{k}_{P}=k_{P}\right) \rightarrow 1.
\end{equation*}

The result implies that for any polynomial $\bar{P}$ with probability
approaching one the lowest order of homogeneous polynomials entering into $%
\bar{P}$ can be determined and also for each coefficient it can be decided
whether it is zero or not with probability approaching 1; each non-zero
coefficient can be consistently estimated.

\subsection{Adaptively estimated asymptotic Wald statistic}

The case of one restriction is given by the following Lemma.

With $q=1$ represent each component $\left\{ G\left( \theta \right) \right\}
_{i}$ of $G\left( \theta \right) $ as a polynomial of form $\left( \ref%
{homog1}\right) $ and consider the corresponding $k_{G_{i}}$ defined in $%
\left( \ref{k-bar1}\right) $ and the corresponding estimator, $\hat{k}%
_{G_{i}}.$ Then define $\hat{k}=\min \left\{ \hat{k}_{G_{i}}\right\} $ and
the corresponding vector $\bar{G}_{\hat{k}}(y)$ with components given by the
homogeneous polynomials of order $\hat{k}$ (some could be zero).

\textbf{Lemma 6.2.} \textit{Under the Assumptions of Theorem 4.1 if (a) }$%
\hat{k}=0,$\textit{\ then with probability approaching 1 as }$T\rightarrow
\infty $\textit{\ there is no singularity, the distribution of the
asymptotic statistic is standard and under the normality assumption is }$%
\chi _{1}^{2};$\textit{\ if (b) }$\hat{k}>0$\textit{\ the estimated
asymptotic statistic is }%
\begin{equation*}
\hat{W}_{T}=\frac{1}{\left( \hat{k}+1\right) ^{2}}\frac{(Z^{\prime }\bar{G}_{%
\hat{k}}\left( Z\right) ^{\prime })^{2}}{\bar{G}_{\hat{k}}\left( Z\right) 
\bar{G}_{\hat{k}}\left( Z\right) ^{\prime }};
\end{equation*}%
and \textit{its distribution converges to the non-standard asymptotic
distribution for the Wald statistic at }$\bar{\theta}$\textit{\ as }$%
T\rightarrow \infty .$

The next theorem considers the general case of the Wald test for several
restrictions. Denote by $\hat{k}_{\det }$\ the estimator of $k_{P}$ applied
to $\bar{P}_{\bar{\theta}}(\theta )$ in $\left( \ref{genericpolyn}\right) $
that represents the polynomial $\det [G\left( \theta \right) G\left( \theta
\right) ^{\prime }];$ as $T\rightarrow \infty $ $\Pr \left( \hat{k}_{\det
}=k_{\det }\right) \rightarrow 1.$ Set $A=I$ then for every $q\times q$
submatrix $\hat{G}_{l}(y)$ of $\hat{G}(y)$ the estimator of $k_{\det ,l}$
for the corresponding determinant polynomial as $T\rightarrow \infty $
equals the true $k_{\det ,l}$ with probability approaching 1, and so does
then the estimated value of $\bar{a}$, as well as the estimated $\alpha _{i}$
defined in proof of Lemma 4.3. It follows that thus one can determine with
probability approaching 1 whether the CLDR property holds and if it does
estimate the matrix $\Lambda _{\bar{\theta}}$ with probability approaching
1. Then the corresponding consistent estimator, $\hat{G}_{\bar{\theta}%
}\left( y\right) $\ of $\bar{G}_{\bar{\theta}}\left( y\right) $ is obtained
for $A=I.$ A consistent estimator of the corresponding polynomial, $\bar{G}_{%
\bar{\theta}}\left( Ay\right) $ for $A=V^{-\frac{1}{2}}$ can be obtained by
substituting $\tilde{y}=\hat{V}^{-\frac{1}{2}}y$ into the estimator $\hat{G}%
_{\bar{\theta}}\left( y\right) $ to obtain $\tilde{G}_{\bar{\theta}}(\tilde{y%
}).$

\textbf{Theorem 6.1.} \textit{Under the Assumptions of Theorem 4.1 if (a)
the corresponding estimated }$\hat{k}_{\det }=0,$\textit{\ then with
probability approaching 1 as }$T\rightarrow \infty $\textit{\ the
distribution of the asymptotic statistic is standard and under normality is }%
$\chi _{q}^{2};$\textit{\ if (b) for }$A=I,$\textit{\ }$\hat{k}_{\det }\neq
0 $\textit{\ but the estimated }$\hat{G}_{\bar{\theta}}\left( y\right) $%
\textit{\ has deficient rank as }$T\rightarrow \infty $\textit{, then the
statistic diverges to infinity; if (c) with }$A=I,$\textit{\ }$\hat{k}_{\det
}\geq 1$\textit{\ and the estimated }$\bar{G}_{\bar{\theta}}\left( y\right) $%
\textit{\ satisfies the CLDR property with the estimated matrix }$\hat{%
\Lambda}_{\bar{\theta}},$\textit{\ then the limit distribution is
consistently estimated by }%
\begin{equation*}
Z^{\prime }[\tilde{G}_{\bar{\theta}}(Z)]^{\prime }\hat{\Lambda}_{\bar{\theta}%
}\left\{ \left[ \tilde{G}_{\bar{\theta}}(Z)\right] [\tilde{G}_{\bar{\theta}%
}(Z)]^{\prime }\right\} ^{-1}\hat{\Lambda}_{\bar{\theta}}[\tilde{G}_{\bar{%
\theta}}(Z)]Z.
\end{equation*}

\subsection{Conservative tests with adaptively estimated bounds}

From Theorem 5.1 it follows that if the CLDR property holds the bound is
provided by 
\begin{equation*}
\frac{1}{(1+\alpha )^{2}}\left\Vert Z^{\prime }Z\right\Vert ,\text{ with }%
\left\Vert Z^{\prime }Z\right\Vert \text{ distributed as }\chi _{p}^{2}\text{
under normality.}
\end{equation*}%
For one restriction CLDR always holds and $\hat{a}=\hat{k}$ as defined in
Lemma 6.2 provides $\alpha $ with probability approaching 1.

For several restrictions estimate $\hat{k}_{\det },$ $\hat{G}_{\bar{\theta}%
}(y)$ as defined in Theorem 6.1. and then if CLDR property holds it is
sufficient to define the estimate of $\alpha $ as the smallest diagonal
element of $\hat{\Lambda}_{\bar{\theta}}$ as in Theorem 6.1. This estimator
will equal the true $\alpha $ with probability approaching 1.

Use the bound $\frac{1}{\left( \hat{a}+1\right) ^{2}}\chi _{p}^{2}.$

\section{Appendix: Proofs}

\begin{proof}[Proof of lemma 4.1]
We first note that for the submatrix $\tilde{G}_{l}(y)$\ defined by (\ref%
{deg}),%
\begin{equation*}
\bar{G}_{l}(y)=\lim_{\lambda =+\infty }diag(\lambda ^{\alpha _{i}})S\tilde{G}%
_{l}(y/\lambda )
\end{equation*}%
is a submatrix of $\bar{G}\left( y\right) .$ Then%
\begin{equation*}
\lim_{\lambda \rightarrow +\infty }\lambda ^{\sum_{i=1}^{q}\alpha _{i}}\det
(S)\det (\tilde{G}_{l}(y/\lambda ))
\end{equation*}%
exists and since $S$ is nonsingular 
\begin{equation*}
\lim_{\lambda \rightarrow +\infty }\lambda ^{\sum_{i=1}^{q}\alpha _{i}}\det (%
\tilde{G}_{l}(y/\lambda ))<\infty .
\end{equation*}%
Then this is%
\begin{equation*}
\lim_{\lambda \rightarrow +\infty }\left( \lambda ^{\sum_{i=1}^{q}\alpha
_{i}-\bar{a}}\right) \lambda ^{\bar{a}}\det (\tilde{G}_{l}(y/\lambda ))
\end{equation*}%
and it follows that $\sum_{i=1}^{q}\alpha _{i}-\bar{a}\leq 0.$ If $%
\sum_{i=1}^{q}\alpha _{i}-\bar{a}=0$ then $\bar{G}_{l}\left( y\right) $ is
full rank and so is $\bar{G}\left( y\right) .$ If $\bar{G}\left( y\right) $
is full rank then there is a square submatrix $\bar{G}_{l}\left( y\right) $
of full rank, for the corresponding submatrix in $SG\left( y\right) $%
\begin{equation*}
\lim_{\lambda \rightarrow +\infty }\lambda ^{\sum_{i=1}^{q}\alpha _{i}}\det (%
\tilde{G}_{l}(y/\lambda ))>0
\end{equation*}%
and $\sum_{i=1}^{q}\alpha _{i}-\bar{a}\geq 0$. The equality follows.
\end{proof}

\begin{proof}[Proof of Lemma 4.2]
By the property $\left( \ref{homog}\right) $ some $a$ for which $\left( \ref%
{deg}\right) $ holds exists and by the Lemma 4.1 either CLDR holds or $a$ is
such that $\sum\limits_{i=1}^{q}\alpha _{i}\neq \bar{\alpha}.$ Then by the
condition on $a^{\prime }$ if CLDR does not hold then for some $i$ we have $%
\alpha _{i}^{\prime }>\alpha _{i},$ or $\alpha _{i}^{\prime }<\alpha _{i}.$
Then for any $\left\{ i,j\right\} $ matrix entry in the matrix $SG(y/\lambda
),$ given by $S_{i\cdot }^{\prime }G_{\cdot j}(y/\lambda )$ (where for a
matrix $A,$ $A_{i\cdot }$ denotes the $i$th row and $A_{\cdot j}$ - the $j$%
th column)%
\begin{equation*}
\lambda ^{\alpha _{i}^{\prime }}S_{i\cdot }^{\prime }G_{\cdot j}(y/\lambda
)=\lambda ^{\alpha _{i}^{\prime }-\alpha _{i}}\lambda ^{\alpha
_{i}}S_{i\cdot }^{\prime }G_{\cdot j}(y/\lambda )
\end{equation*}

In the first case $\alpha _{i}^{\prime }>\alpha _{i}$, this matrix entry
diverges to infinity. In the second case $\alpha _{i}^{\prime }<\alpha _{i}$
and as $\lambda \rightarrow \infty $ this matrix entry converges to zero for
every $j,$ thus the limit matrix $\lim_{\lambda \rightarrow +\infty }\left[
diag(\lambda ^{\alpha _{i}^{\prime }})SG(y/\lambda )\right] $ has deficient
rank.
\end{proof}

\begin{proof}[Proof of Lemma 4.3]
Start with a $q_{v}\times p$ matrix of polynomials $G^{v}\left( y\right) .$

For each matrix element, $\left\{ G^{v}\left( y\right) \right\} _{ij},$which
is a polynomial, define the lowest order of homogeneous polynomial, $\bar{k}%
_{\left\{ G^{v}\left( y\right) \right\} _{ij}}.$ Then select $\bar{k}_{v}=%
\underset{i,j}{\min }\{\bar{k}_{\left\{ \ G^{v}\left( \ y\right) \ \right\}
\ _{ij}}\}.$ Consider a polynomial matrix, $\tilde{G}_{v}\left( y\right) $
such that 
\begin{equation*}
\left\{ \tilde{G}_{v}\left( y\right) \right\} =\left\{ 
\begin{array}{cc}
\left\{ G^{v}\left( y\right) \right\} _{ij,\bar{k}_{v}} & \text{when this
polynomial is non-zero} \\ 
0 & \text{otherwise.}%
\end{array}%
\right. 
\end{equation*}%
In other words, the $ij$ element of $\tilde{G}_{v}\left( y\right) $ is
either a non-zero polynomial of order $\bar{k}_{v},$ that entered into $%
\left\{ G^{v}\left( y\right) \right\} _{ij},$ or zero. Next, consider all
square submatrices $r_{v}\times r_{v},$ $r_{v}\leq q_{v}$ of $\tilde{G}%
_{v}\left( y\right) ,$ for at least one of those determinant is non-zero;
select the largest $\bar{r}_{v}$ with the property that some submatrix of
this dimension has a non-zero determinant, and (i) either $\bar{r}_{v}=q_{v},
$ or (ii) determinant of any submatrix with $q_{v}\geq r_{v}>\bar{r}_{v}$ is
zero.

In case (i) define $S^{v}=I_{q_{v}}.$ In case (ii) construct a non-singular
matrix $S^{v},$ such that for some $\bar{r}_{v}\times p$ full row rank
matrix of polynomials, $\bar{G}^{v}\left( y\right) ,$ 
\begin{equation*}
S^{v}\tilde{G}_{v}\left( y\right) =\left[ 
\begin{array}{c}
\bar{G}^{v}(y) \\ 
0%
\end{array}%
\right] .
\end{equation*}%
Such a matrix always exists. Then $S^{v}G^{v}\left( y\right) $ has the
representation%
\begin{equation*}
\left[ 
\begin{array}{c}
\bar{G}^{v}(y)+N^{v}\left( y\right) \\ 
G^{v+1}\left( y\right)%
\end{array}%
\right]
\end{equation*}%
where if $N^{v}\left( y\right) $ is non-zero the polynomial entries in the
matrix $N^{v}\left( y\right) $ have homogeneous polynomial terms of order no
less than $\bar{k}_{v};$ and the non-zero $\left( q_{v}-\bar{r}_{v}\right)
\times p$ matrix $G^{v+1}\left( y\right) $ has polynomial terms only of
order $\geq \bar{k}_{v}+1.$

Consider now the original matrix $G\left( y\right) ,$ denote it $G^{1}\left(
y\right) $ with $q_{1}=q$ and employ the construction recursively until for
some $m$ it ends: $\Sigma _{v=1}^{m}\bar{r}_{v}=q.$

Denote by $\bar{S}^{v}$ the matrix $\left[ 
\begin{array}{c}
I_{\bar{r}_{1}+...+\bar{r}_{v-1}} \\ 
S^{v}%
\end{array}%
\right] $ and define $S=\bar{S}^{m}...\bar{S}^{1}.$ Set $a=\left( \alpha
_{1},...,\alpha _{q}\right) =(\bar{k}_{1},...,\bar{k}_{1},...,\bar{k}%
_{m},...,\bar{k}_{m}),$ where each $\bar{k}_{v}$ enters $\bar{r}_{v}$ times.
Then for this $a$ and $S$ 
\begin{equation*}
\lim_{\lambda \rightarrow \infty }\left[ diag(\lambda ^{\alpha
_{i}})SG(y/\lambda )\right]
\end{equation*}%
is a finite matrix $\bar{G}\left( y\right) =\left[ 
\begin{array}{c}
\bar{G}^{1}\left( y\right) \\ 
\vdots \\ 
\bar{G}^{m}\left( y\right)%
\end{array}%
\right] ;$ if $\sum\limits_{i=1}^{q}\alpha _{i}=\bar{\alpha},$ then CLDR
property holds, if $\sum\limits_{i=1}^{q}\alpha _{i}<\bar{\alpha}$ then the
limit matrix has deficient rank.
\end{proof}

\begin{proof}[Proof of Theorem 4.1]
Consider $y_{T}^{\ast }=\lambda _{T}y_{T}$ and the quadratic form similar to 
$\left( \ref{quad}\right) $ 
\begin{equation*}
W(y_{T}^{\ast },g_{\bar{\theta}},\lambda _{T},A\hat{V}_{T}A^{\prime
})=\lambda _{T}^{2}g_{\bar{\theta}}^{\prime }(y_{T}^{\ast }/\lambda _{T})[G_{%
\bar{\theta}}(y_{T}^{\ast }/\lambda _{T})A\hat{V}_{T}A^{\prime }G_{\bar{%
\theta}}^{\prime }(y_{T}^{\ast }/\lambda _{T})]^{-1}g_{\bar{\theta}%
}(y_{T}^{\ast }/\lambda _{T}).
\end{equation*}%
From Assumption 2.3 if $\lambda =\lambda _{T}$ and $\theta =\hat{\theta}_{T}$
then the probability limit of corresponding $V^{-\frac{1}{2}%
}A^{-1}y_{T}^{\ast }$ is $Z$ with distribution $Q\left( \bar{\theta}\right) ;
$ from Assumption 2.4 $\hat{V}_{T}=V+o_{p}(1).$ From $\left( \ref{decomp2}%
\right) $ and convergence it follows that 
\begin{eqnarray}
diag(\lambda _{T}^{\alpha _{i}})SG_{\bar{\theta}}(y_{T}^{\ast }/\lambda
_{T}) &=&\bar{G}_{\bar{\theta}}(y_{T}^{\ast })+O_{p}(1/\lambda _{T}); \\
diag(\lambda _{T}^{\alpha _{i}})S\lambda g_{\bar{\theta}}(y_{T}^{\ast
}/\lambda _{T}) &=&\bar{g}_{\bar{\theta}}(y_{T}^{\ast })+O_{p}(1/\lambda
_{T}).  \notag
\end{eqnarray}%
Then $W(y_{T}^{\ast },g_{\bar{\theta}},\lambda _{T},A\hat{V}_{T}A^{\prime })=
$%
\begin{eqnarray*}
&&\lbrack \bar{g}_{\bar{\theta}}(y_{T}^{\ast })+O_{p}(1/\lambda
_{T})]^{\prime }\left\{ \left[ \bar{G}_{\bar{\theta}}(y_{T}^{\ast
})+O_{p}(1/\lambda _{T})\right] A[V+o_{p}(1)]A^{\prime }\left[ \bar{G}_{\bar{%
\theta}}(y_{T}^{\ast })+O_{p}(1/\lambda _{T})\right] ^{\prime }\right\} 
^{\substack{ - \\ -1}} \\
&&\times \lbrack \bar{g}_{\bar{\theta}}(y_{T}^{\ast })+O_{p}(1/\lambda
_{T})].
\end{eqnarray*}

(a) If CLDR holds then $W_{T}$ by continuity of the determinants of
polynomials and polynomial matrices converges to 
\begin{equation*}
\lbrack \bar{g}_{\bar{\theta}}(Z)]^{\prime }\left\{ [\bar{G}_{\bar{\theta}%
}(Z)]AVA^{\prime }[\bar{G}_{\bar{\theta}}(Z)]^{\prime }\right\} ^{-1}[\bar{g}%
_{\bar{\theta}}(Z)];
\end{equation*}%
substituting the reparametrized functions for $A=V^{-\frac{1}{2}}$ we get
the result.

(b) Follows by continuity of the determinants of polynomial matrices and $%
\left( \ref{blowup}\right) .$
\end{proof}

\begin{proof}[Proof of Theorem 5.1]
Consider the asymptotically equivalent statistic:%
\begin{eqnarray*}
&&Z^{\prime }\bar{G}(Z)^{\prime }\Lambda \left[ \bar{G}(Z)\bar{G}(Z)^{\prime
}\right] ^{-1}\Lambda \bar{G}(Z)Z \\
&=&Z^{\prime }\bar{G}(Z)^{\prime }\left( \bar{G}(Z)\bar{G}(Z)^{\prime
}\right) ^{-\frac{1}{2}} \\
&&\times \left[ \left( \bar{G}(Z)\bar{G}(Z)^{\prime }\right) ^{\frac{1}{2}%
}\Lambda \left( \bar{G}(Z)\bar{G}(Z)^{\prime }\right) ^{-\frac{1}{2}}\left( 
\bar{G}(Z)\bar{G}(Z)^{\prime }\right) ^{-\frac{1}{2}}\Lambda \left( \bar{G}%
(Z)\bar{G}(Z)^{\prime }\right) ^{\frac{1}{2}}\right] \\
&&\times \left( \bar{G}(Z)\bar{G}(Z)^{\prime }\right) ^{-\frac{1}{2}}\bar{G}%
(Z)Z \\
&\leq &\left\Vert \left( \bar{G}(Z)\bar{G}(Z)^{\prime }\right) ^{-\frac{1}{2}%
}\bar{G}(Z)Z\right\Vert ^{2}\left\Vert \left( \bar{G}(Z)\bar{G}(Z)^{\prime
}\right) ^{\frac{1}{2}}\Lambda \left( \bar{G}(Z)\bar{G}(Z)^{\prime }\right)
^{-\frac{1}{2}}\right\Vert \\
&&\times \left\Vert \left( \bar{G}(Z)\bar{G}(Z)^{\prime }\right) ^{-\frac{1}{%
2}}\Lambda \left( \bar{G}(Z)\bar{G}(Z)^{\prime }\right) ^{\frac{1}{2}%
}\right\Vert \\
&\leq &\left\Vert \Lambda \right\Vert ^{2}\left\Vert Z\right\Vert ^{2}\sim 
\frac{1}{\left( 1+i_{0}\right) ^{2}}\chi _{p}^{2},
\end{eqnarray*}%
since for similar matrices the eigenvalues are the same, so eigenvalues of $%
\left( \bar{G}(Z)\bar{G}(Z)^{\prime }\right) ^{\frac{1}{2}}\Lambda \left( 
\bar{G}(Z)\bar{G}(Z)^{\prime }\right) ^{-\frac{1}{2}}$ are the same as for $%
\Lambda $ regardless of $Z$ and the norm is given by the largest eigenvalue,
and finally, 
\begin{equation*}
\left\Vert \left( \bar{G}(Z)\bar{G}(Z)^{\prime }\right) ^{-\frac{1}{2}}\bar{G%
}(Z)Z\right\Vert ^{2}=\left( Z^{\prime }\bar{G}(Z)^{\prime }\left( \bar{G}(Z)%
\bar{G}(Z)^{\prime }\right) ^{-1}\bar{G}(Z)Z\right)
\end{equation*}%
where for every value of $Z$ the corresponding constant matrix $\bar{G}%
^{\prime }(Z)\left( \bar{G}(Z)\bar{G}^{\prime }(Z)\right) ^{-1}\bar{G}(Z)$
is a projection and thus its norm is always bounded by 1.
\end{proof}

\begin{proof}[Proof of Lemma 5.1]
Express the p.d.f. of $\chi _{p_{1}}^{2}/\alpha _{1}:$%
\begin{equation*}
p.d.f._{\chi _{p_{1}}^{2}/\alpha _{1}}(y)=\frac{\alpha _{1}}{%
2^{p_{1}/2}\Gamma \left( \frac{p_{1}}{2}\right) }\exp (-\alpha
_{1}y/2)\left( \alpha _{1}y\right) ^{\frac{p_{1}}{2}-1},
\end{equation*}%
and similarly for $\chi _{p_{2}}^{2}/\alpha _{2}.$ The ratio $\frac{%
p.d.f._{\chi _{p_{1}}^{2}/\alpha _{1}}(y)}{p.d.f._{\chi _{p_{2}}^{2}/\alpha
_{2}}(y)}$ is%
\begin{equation*}
2^{\frac{p_{2}-p_{1}}{2}}\left( \Gamma \left( \frac{p_{2}}{2}\right) /\Gamma
\left( \frac{p_{1}}{2}\right) \right) \frac{\alpha _{2}^{\frac{p_{2}}{2}-1}}{%
\alpha _{1}^{\frac{p_{1}}{2}-1}}y^{\frac{p_{1}-p_{2}}{2}}\exp \left( \frac{y%
}{2}(\alpha _{2}-\alpha _{1})\right) .
\end{equation*}%
Since $\alpha _{2}>\alpha _{1}$ for large enough $y$ this expression is
larger than 1.
\end{proof}

\begin{proof}[Proof of Lemma 6.1]
First consider $\hat{P}\left( i_{1},...,i_{p},\bar{\theta}\right) $ as
defined in $\left( \ref{estimator}\right) .$ By polynomial structure and the
convergence rate in Assumption 2.3. $\hat{P}\left( i_{1},...,i_{p},\bar{%
\theta}\right) =P\left( i_{1},...,i_{p},\bar{\theta}\right) +O_{p}\left(
\lambda ^{-1}\right) .$ Two consequence are (a) from Assumption 2.4 then $%
\hat{P}\left( i_{1},...,i_{p},\bar{\theta},V\right) -\bar{P}\left(
i_{1},...,i_{p},\bar{\theta},V\right) \rightarrow _{p}0;$ (b) when $P\left(
i_{1},...,i_{p},\bar{\theta}\right) =0,$ $\Pr \left( \hat{P}\left(
i_{1},...,i_{p},\bar{\theta}\right) =0\right) \rightarrow 1$ by construction 
$\left( \ref{estimator}\right) .$ Since $\hat{k}_{P}-1$ can be defined as
the highest order of polynomial with $\hat{P}\left( i_{1},...,i_{p},\bar{%
\theta}\right) =0$ it follows $\Pr \left( \hat{k}_{P}=k_{P}\right)
\rightarrow 1;$ note that $\hat{k}_{P}$ for a polynomial constructed for $%
A=I $ is the same as for any non-singular $A.$
\end{proof}

\begin{proof}[Proof of Lemma 6.2]
Apply Lemma 6.1 to each of the estimated polynomials to determine with
probability approaching 1 the lowest order $k_{P}$ of the non-zero
homogeneous polynomial and to obtain the consistent estimators of the
polynomial vector functions, $\bar{G}\left( \cdot \right) .$ Substituting
the limit $Z$ provides the consistent estimator of the asymptotic
distribution.
\end{proof}

\begin{proof}[Proof of Theorem 6.1]
The proof follows by application of Lemma 1 to each polynomial that is
estimated.
\end{proof}

\end{document}